\documentclass[a4paper, 11pt]{article}
\usepackage[english]{babel}
\usepackage[T1]{fontenc}
\usepackage{amsmath}
\usepackage{amsthm}
\usepackage{amssymb}
\usepackage{amsfonts}
\usepackage{graphicx}
\usepackage{caption}
\usepackage{xcolor}
\usepackage{subfig}

\graphicspath{{img/}}

\begin{document}
\title{Overcoming slowly decaying Kolmogorov $n$-width by transport maps: application to model order reduction of fluid dynamics and fluid--structure interaction problems}

\author{Monica Nonino, Francesco Ballarin, Gianluigi Rozza\\ Mathematics Area, mathLab, SISSA, Trieste, Italy \and Yvon Maday\\ Sorbonne Universit\'e, CNRS, Universit\'e de Paris, \\Laboratoire Jacques-Louis Lions, \\F-75005 Paris, France and Institut Universitaire de France}

\maketitle
\abstract{In this work we focus on reduced order modelling for problems for which the resulting reduced basis spaces show a slow decay of the Kolmogorov $n$-width, or, in practical calculations, its computational surrogate given by the magnitude of the eigenvalues returned by a proper orthogonal decomposition on the solution manifold. In particular, we employ an additional preprocessing during the offline phase of the reduced basis method, in order
to obtain smaller reduced basis spaces. Such preprocessing is based on the composition of the snapshots with a transport map, that is a family of smooth and invertible mappings that map the physical domain of the problem into itself.
Two test cases are considered: a fluid moving in a domain with deforming walls, and a fluid past a rotating cylinder.
Comparison between the results of the novel offline stage and the standard one is presented.}\\

\noindent\small{\textbf{keywords:} reduced order modeling; transport dominated problems; fluid-structure interaction problem; flow past a rotating cylinder; proper orthogonal decomposition.}
\section{Introduction}
The reduced basis method \cite{HesthavenRozzaStamm,Haasdonk2017,RozzaHuynhPatera1,BertagnaVeneziani,Lassila,Haasdonk,Nguyen}
is a powerful tool when requiring fast simulations of parametrized partial differential equations (PDEs): its efficiency relies on the possibility to construct an approximation of the solution for any value of the parameter in the span of a few
basis functions, which are computed in the (expensive) offline phase. 
Despite its capability has been acclaimed in a large variety of situations, model reduction of advection dominated (or even hyperbolic) problems is still a challenging task \cite{Ohlberger2016,Greif2019,Ehrlacher2019}. It has therefore become clear that a modification of the way the reduced basis method works is necessary, especially in order to be able to obtain a small basis set also in these more challenging situations.

Assume that $\Omega$ is the physical domain of the problem of interest. Let $\mathcal{P} \subset \mathbb{R}^p$, $p \in \mathbb{N}$, be a compact set, and denote by $\mu\in\mathcal{P}$ the parameter. Furthermore, let $t \in [0, T]$ be the time, for some $T > 0$. In the following the symbol $\eta$ will either stand for $\mu$ in parametrized stationary problems, $t$ in non-parametric unsteady problems, or the pair $(t, \mu)$ in parametrized unsteady problems. For any value of $\eta \in \mathcal{E}$, we seek the solution $\mathbf{z}(\eta)\colon\Omega \to \mathbb{R}^Z$, $Z \in \mathbb{N}$, to the following PDE:
\begin{equation}
N(\mathbf{z}(\eta); \eta) = 0
\label{eq:probN}
\end{equation}
where $N$ is a nonlinear operator working on functions defined over $\Omega$, and $\mathcal{E} := \mathcal{P}$, $\mathcal{E} := [0, T]$ or $\mathcal{E} := [0, T] \times \mathcal{P}$, respectively, in the three aforementioned cases. Specific examples of problems of interests in computational fluid dynamics and fluid-structure interaction will be introduced throughout the paper.

Let $z(\eta)$ be a component of $\mathbf{z}(\eta)$, and $\mathcal{M}_z$ be the solution manifold, embedded in some normed linear space $(X_z, \lvert\lvert\cdot\rvert\rvert_{X_z})$, defined as:
\begin{equation*}
\mathcal{M}_z = \{z(\eta), \quad \eta\in\mathcal{E} \}.
\end{equation*}
The choice of a component-wise solution manifold (rather than just one solution manifold associated to the vector $\mathbf{z}(\eta)$) will be motivated in section 2.
One fundamental assumption of the reduced basis method is that $\mathcal{M}_z$ can be approximated in an accurate way by a sequence of finite dimensional spaces: any element of $\mathcal{M}_z$ can be recovered using a linear combination
of solutions of \eqref{eq:probN}, which are computed only once and for all. The mathematical entity that incorporates this concept is the Kolmogorov $n$-width of $\mathcal{M}_z$, defined as
\begin{equation}
D_n(\mathcal{M}_z, X_z) = \inf\limits_{E_n\subset X_z}\sup\limits_{f\in\mathcal{M}_z}\inf\limits_{g\in E_n} ||f-g||_{X_z},
\end{equation}
where $E_n$ is any linear subspace of $X_z$ of dimension $n$.
The Kolmogorov $n$-width $D_n$ tells us the entity of the error that we commit by approximating any element $f$ of $\mathcal{M}_z$ with an element $g$ of a linear space $E_n$.
Therefore the faster $D_n$ decays as we let $n$ grow, the greater possibility we have to build a good linear approximation space of $\mathcal{M}_z$ of low dimension.
In the majority of the problems there is no explicit analytic formula for $D_n$, yet there are some situations where we can compute good bounds on the Kolmogorov $n$-width\cite{Melenk, CohenDeVore}. 
In general, since it is very difficult to provide such bounds on $D_n$, we can only hope that the $n$-width of the solution manifold is small. 
A heuristic way to check that this hypothesis is satisfied by the problem of interest is to run a proper orthogonal decomposition (POD) on a set of snapshots, and check the rate of decay of the eigenvalues $\{\lambda_i\}_i$ returned by the POD:
if the $\{\lambda_i\}_i$ decay fast, then we can expect to be able to build a good low dimensional linear approximation space for $\mathcal{M}_z$.
This assumption often fails in transport-dominated problems, which show a very slow decay of the eigenvalues of the POD, and thus the inability of the reduced basis method to reconstruct any element of $\mathcal{M}_z$ by using a small number of basis functions.

A growing number of works which focus on constructing alternative (nonlinear) model order reduction techniques, that can be effectively applied to transport-dominated problems, has appeared in recent years.
Abgrall et al. \cite{Abgrall2016} show an $L^1$ norm minimization technique, to be used for the approximation of nonlinear hyperbolic equations, without anyway curing the problem of high dimensional solution manifolds. 
Approximated Lax pairs have been employed for model reduction of nonlinear problems arising in cardiac electrophysiology\cite{Gerbeau2014,Gerbeau2015} . Gerbeau et al \cite{Gerbeau2014} construct a reduced basis space
at each time step: the reduced basis are the modes of a Schroedinger operator where the potential is the solution at the previous time step. The drawback of this procedure is the increase of the number of basis functions as the
approximation accuracy requirement is increased.
Nonlinear model reduction techniques in metric spaces are proposed by Ehrlacher et al. \cite{Ehrlacher2019}, relying either on a tangent principal component analysis or barycentric greedy algorithm. 
Adaptivity is also employed to overcome difficulties arising in model reduction of transport dominated flows, see Carlberg \cite{Carlberg} for the use of localized reduced bases or Peherstorfer \cite{Peherstorfer}.
the procedure to update the basis set.
Several other different techniques rely on a composition with a suitably defined transport map (on which we will focus on in the rest of the work). Such map may be obtained as the solution of a Monge-Kantorovich optimal 
mass transport problem \cite{IolloLombardi,BeIoRi2018}, or provided analytically in simple cases as the one-dimensional problem considered by Cagniart et al. \cite{CagniartMadayStamm}. 
More complicated configurations, possibly including shocks, can still be handled relying on transport maps based on more advanced shape parametrization maps \cite{CagniartCrisovanMadayAbgrall}. 
Extension to multiple transport phenomena are also possible, most notably by means of the shifted POD \cite{Shifted_POD}, as well as (possibly interacting) shocks by recent extensions of the transported snapshots interpolation 
\cite{Welper2017,Welper2017b,Welper2019} and transported snapshots model reduction methods \cite{Nair2019}. 
While in simple cases the underlying transport maps are provided by the user, registration techniques can be employed to automate their selection in more complex geometrical configurations \cite{Taddei}. 
Application of the use of transport maps for model reduction based on an embedded high fidelity method is shown by Karatzas et al. \cite{KaratzasBallarinRozza2019}.
An alternative is the freezing method \cite{Freezing_method,Rowley_2003,BeynThummler}, in which the key tool is the identification of a Lie group acting on a frozen solution component. 
However, to properly develop and analyze the
proposed methodology for complex problems in fluid dynamics, very involved mathematical tools and settings are needed which as of now hinders its applicability in a broader setting. 
This consideration, together with previous observations on other techniques, makes it clear that there is the need for a lighter, simpler and more natural framework. A methodology that satisfies these requirements 
and that is based on the definition of some transport maps has been
introduced and applied to some toy problems by Cagniart \cite{tesidottoratocagniart} and Cagniart et al.\cite{CagniartMadayStamm}. 
The goal of this work is to present an application of model reduction based on transported snapshots for problems in fluid dynamics and fluid-structure interaction, focusing in particular on the comparison between results of the standard offline phase and the novel one which relies on transport. As (especially in fluid dynamics problems) further challenges are present during the online stage, such as keeping into account stabilization \cite{TorloBallarinRozza2018,Azaiez2019} or turbulence modelling \cite{ChaconDelgadoGomezBallarinRozza2017,StabileBallarinZuccarinoRozza2018}, we limit our exposition to the offline stage; future research work will extend the results presented here to the online stage. The article is organized as follows: in Section $2$ we present the nonlinear model reduction approach by transport maps, and we explain how it helps to overcome the problem of slow decay of the Kolmogorov $n$-width. In Section
$3$ we present a first example of application: a time dependent CFD problem, without parameters. In Section $4$ we study a fluid--structure interaction problem, where the solution behaves like a travelling wave, thus
featuring a slow decay of the Kolmogorov $n$-width of the problem; also in this case, the problem is time dependent and non-parametric. In Section $5$, we further extend the CFD test case introduced in Section $3$ by adding a physical parameter, namely the Reynolds number. Conclusions follow in section 6.

\section{Nonlinear model reduction by transport maps}
In this section we summarize the nonlinear approach that we will apply in the forthcoming sections to fluid dynamics and fluid-structure interaction problems. 
We will closely follow the presentation and the notation introduced by Cagniart et al. \cite{CagniartMadayStamm} for a parametrized viscous Burgers equation.
The idea proposed therein is to ``pre-process'' the solution manifold $\mathcal{M}_z$ by a composition with a map in a family of smooth and invertible mappings 
\begin{equation*}
\mathcal{F}_z = \{ F_\eta: \Omega\to\Omega,\ F_\eta \text{ is smooth and invertible},\ \eta \in \mathcal{E} \}.
\end{equation*}
Maps in the family $\mathcal{F}_z$ are parametrized by the same $\eta$ appearing in \eqref{eq:probN}, and are essentially problem-specific. For instance, in the work by Cagniart et al. \cite{CagniartMadayStamm}, a family of translations $F_\eta(x) = x - \eta, x \in \Omega \subset \mathbb{R}$ was employed for the viscous Burgers equation. Specific choices of $\mathcal{F}_z$ will be discussed alongside the problem of interest. We then introduce the ``pre-processed'' solution manifold
\begin{equation*}
\mathcal{M}_{\mathcal{F}_z}=\{z(\eta) \circ F^{-1}_{\eta}, \quad \eta \in \mathcal{E}\}.
\end{equation*}
Assuming that $\mathcal{F}_z$ is carefully chosen, $\mathcal{M}_{\mathcal{F}_z}$ has a smaller Kolmogorov $n$-width than $\mathcal{M}_{z}$. We remark here that one may well choose different families $\mathcal{F}_z$ for different components $z(\eta)$ of the solution $\mathbf{z}(\eta)$; indeed, it is often the case, especially in multiphysics problems, that the qualitative evolution in $\eta$ (which usually suggests the choice of the family $\mathcal{F}_z$ itself) behaves differently for different components of the solution.

The practical realization of this preprocessing procedure is incorporated in the offline phase. Given a discrete training set $\mathcal{E}_{tr}$, we compute each solution component $z(\eta_{tr})$ associated to any $\eta_{tr} \in \mathcal{E}_{tr}$. The discrete approximations of the corresponding standard and preprocessed solution manifolds
\begin{align*}
&\mathcal{M}_z^{tr} = \{z(\eta_{tr}), \quad \eta_{tr}\in\mathcal{E}_{tr} \},\\
&\mathcal{M}_{\mathcal{F}_z}^{tr}=\{z(\eta_{tr}) \circ F^{-1}_{\eta_{tr}}, \quad \eta_{tr} \in \mathcal{E}_{tr}\},
\end{align*}
provide snapshots for a compression by a POD. 
The compression is here applied to both $\mathcal{M}_z^{tr}$ and $\mathcal{M}_{\mathcal{F}_z}^{tr}$ to provide a comparison between the standard offline phase and one with preprocessing, but in practical computations 
one would neglect the compression of $\mathcal{M}_z^{tr}$, as it is understood that $\mathcal{M}_{\mathcal{F}_z}^{tr}$ would result in a POD basis set $\{\Phi_i\}_{i=1}^N$ of lower dimension. 
\\As anticipated in the previous section, we will only focus on the offline part of the reduced order method; efficient evaluation of the online solution in a Galerkin projection setting is a future development of this work; 
nevertheless we provide for completeness some details on how the online phase of this new reduction method should be carried out. Let us assume that an accurate approximation $z^N(t^n)$ of the solution component $z(t^n)$ at time $t^n$ is known,
as a linear combination of our basis functions $\Phi_i$:
\begin{equation*}
z^N(t^n) = \sum_{i=1}^N \alpha_i^n\Phi_i\circ F_{\eta_{tr}^n}.
\end{equation*}
In order to recover $z^N(t^{n+1})$ as an expansion:
\begin{equation*}
z^N(t^{n+1}) = \sum_{i=1}^N \alpha_i^{n+1}\Phi_i\circ F_{\eta_{tr}^{n+1}},
\end{equation*}
the idea proposed by Cagniart et al.\cite{CagniartMadayStamm} is to iterate between the search for the reduced coordinates $\alpha_i^{n+1}$ and the search for a suitable parameter $\eta_{tr}^{n+1}$. 
This procedure is carried out by means of a minimization problem of the $L^2$ norm of the residual evaluated at $\sum_{i=1}^N \alpha_i^{n+1}\Phi_i\circ F_{\eta_{tr}^{n+1}}$.

In the next two sections we will see two applications of this preprocessing procedure: the first problem we are going to study is a computational fluid dynamics (CFD) problem, and the second problem is a fluid--structure interaction problem. 
These two applications are quite different for what concerns the physics behind them (the latter is a coupled multiphysics problem, the former is not), but they both feature a slow decay of the eigenvalues returned by running
a POD on $\mathcal{M}_z^{tr}$.
In the first case is characterized by the change of the direction of propagation of a vortex close to a rotating cylinder, a feature which is difficult to reproduce with a small number of basis functions.
The second problem is a transport dominated problem, where the solution behaves like a wave travelling in the domain. Also in this case the travelling wave would be hard to be reconstructed with a standard model order reduction technique.

\section{A CFD test case}
In this Section we are going to show the results that we obtain with the preprocessing procedure on a CFD test case, namely a fluid past a rotating cylinder. 
The problem described in the following is inspired by the one presented in Cagniart's thesis \cite{tesidottoratocagniart}, although our formulation is slightly different: in the reference, the direction of the fluid velocity at the inlet boundary is changing 
in time, and the cylinder is kept fixed, whereas on the contrary in our problem the inlet direction is constant, but the cylinder is rotating.
\subsection{Problem formulation}
In this problem we have a flow past a rotating cylinder; this situation, with or without the rotation feature, is quite interesting and has been indeed studied in a large number of papers, both in the incompressible case \cite{Rotating_cylinder_Mittal_Behara,High_rotation_rates_Breuer,cinesi,StojkovicSchoenBreuerDurst} and in the compressible regime \cite{TeymourtashSalimipour}. 
It is well known that if the Reynolds number $Re$ is greater than $47$, then a vortex shedding phenomenon occurs \cite{cinesi}. When, in addition, the cylinder is rotating, it might happen that we see 
a change in the direction of propagation of this vortex; this phenomenon is quite complicated and is strictly related to a variety of different physical quantities. In fact, it is related not only to $Re$, but also to the cylinder rotation rate $\alpha$, 
which is defined as:
\begin{equation*}
\alpha = \frac{D\omega}{2U_{\infty}},
\end{equation*}
where $D$ is the diameter of the cylinder, $\omega$ is the angular velocity of a point on the surface of the cylinder and $U_{\infty}$ is the oncoming free stream velocity. For $\alpha < \alpha_L$ vortex shedding occurs,
where $\alpha_L$ is a critical value that is a function of $Re$. We refer to Stojkovi\'c et al. \cite{High_rotation_rates_Breuer} for different values of $\alpha$ related to different values of the Reynolds number.

\begin{figure}
\centering
\includegraphics[scale=0.2]{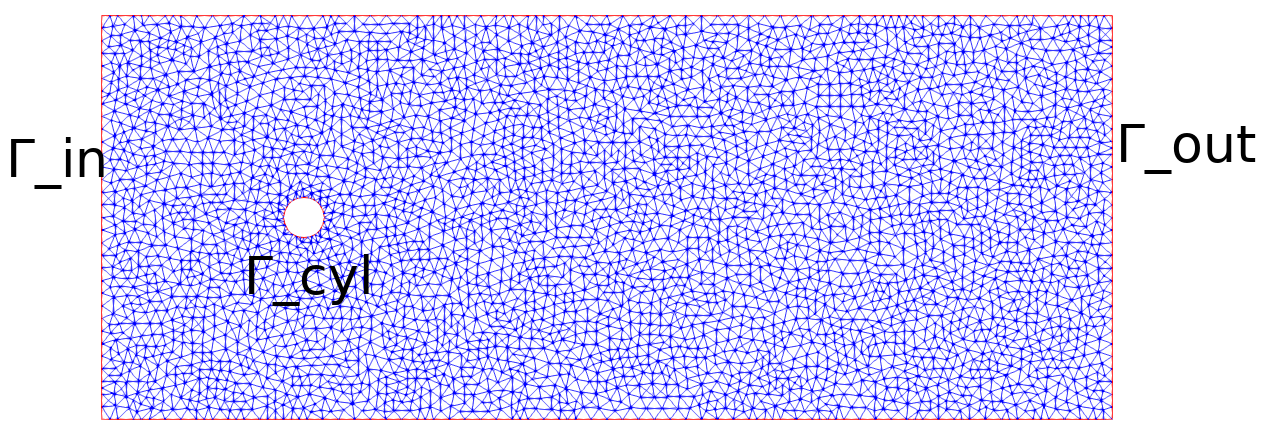}
\caption{Physical domain of the CFD test case.}
\label{CFD mesh}
\end{figure}

To model our problem we use incompressible Navier-Stokes equation. Figure \ref{CFD mesh} shows the physical domain $\Omega$ for our problem.
We do not take into account any physical
parameter. Our problem reads as follows: for any time $t \in [0, T]$ find $\textbf{u}_f(\cdot; t)$ and $p_f(\cdot; t)$ such that:
\begin{equation*}
\begin{cases}
\partial_t\textbf{u}_f + (\textbf{u}_f\cdot\nabla)\textbf{u}_f - \frac{1}{Re}\Delta\textbf{u}_f + \nabla p_f = \textbf{b}_f &\text{in $\Omega \times [0, T]$,}\\
-\text{div}\textbf{u}_f = 0 &\text{in $\Omega \times [0, T]$,}\\
\textbf{u}_f = \textbf{u}_{in} &\text{in $\Gamma_{in}\times [0, T]$,}\\
\textbf{u}_f = \textbf{u}_{tan} &\text{in $\Gamma_{cyl}\times [0, T]$,}\\
p_f \mathbf{n} - \frac{1}{Re} \nabla\textbf{u}_f \cdot \mathbf{n} = \mathbf{0} &\text{in $\left(\Gamma_{top} \cup \Gamma_{bottom} \cup \Gamma_{out}\right) \times [0, T]$.}
\end{cases}
\end{equation*}
We impose homogeneous Neumann conditions on top and bottom walls $\Gamma_{top}$ and $\Gamma_{bottom}$, and also on the right boundary $\Gamma_{out}$, as all these boundaries are considered as outlets due to the fact that the vortex sheet rotates alongside with the cylinder. We impose a 
Dirichlet condition $\textbf{u}_f = \textbf{u}_{in}$ on $\Gamma_{in}$, where $\textbf{u}_{in}$ is a fixed horizontal inflow tabulated in Table \ref{CFD parameters}. 

Furthermore, $\textbf{u}_{tan}$ is the tangential velocity at the surface of the rotating cylinder.
At the beginning of the simulation the cylinder is not moving, and it stays still until a fully developed vortex shedding phenomenon is reached; after that the cylinder starts to rotate counterclockwise, first with
a constant acceleration $\beta$, then it keeps on rotating with a constant angular velocity. 
Thus, denoting by $\omega$ be the angular velocity of the cylinder, we set
\begin{equation}
\omega =
\begin{cases}
\omega_0 = 0 &\text{for $t \in [0, t_1]$,}\\
\omega_t = \omega_0 + \beta (t-t_1) &\text{for $t \in [t_1, t_2]$,}\\
\omega_f = \omega_{t_2} &\text{for $t \in [t_2, T]$.} 
\end{cases}
\end{equation}
Once we have the angular velocity, we can compute $\textbf{u}_{tan}$ thanks to the relation:
\begin{equation*}
\lvert\lvert \textbf{u}_{tan}(t)\rvert\rvert = \omega_t r,
\end{equation*}
where $r$ is the radius of the cylinder, and assuming $\textbf{u}_{tan}(t)$ to be tangent to $\Gamma_{cyl}$.\\
As far as the problem formulation concerns, let us mention the fact that, in preparation of a forthcoming online stage, we adopted a supremizer enrichment technique, meaning that we decided to enrich the space 
of the fluid velocity snapshots $\textbf{u}_f$ with some supremizers snapshots $s_u$. The supremizer snapshots $s_u$ are needed at the reduced order level 
in order to have a more stable approximation of the fluid pressure. 
For further details on the formulation of the supremizer in the POD framework of parametrized fluid flows we refer to Ballarin et al.\cite{Supremizer}.
\begin{table}
\begin{center}
\begin{tabular}{|l|c||l|c|}
\hline
Data & Value & Data & Value\\
\hline
$r$ & $2$ cm & $\textbf{u}_{in}$ & $(1, 0)$\\ 
\hline
$Re$ & $100$ & $\Delta t$ & $0.25$\\ 
\hline 
$b_f$ & $0$ & $\beta$ & $0.025$ rad/s$^2$\\
\hline
$t_1$ & $75$ s & $t_2$ & $95$ s\\
\hline
$T$ & $145$ s & & \\
\hline
\end{tabular}
\end{center}
\caption{Problem data for the CFD test case.}
\label{CFD parameters}
\end{table}

\begin{figure}
\centering
\includegraphics[scale=0.5]{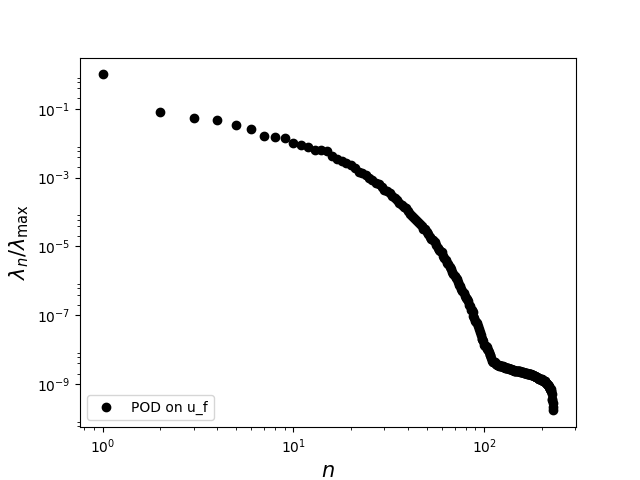}
\caption{Decay of the eigenvalues for the fluid velocity in the simulation with a rotating cylinder.}
\label{POD_Karman}
\end{figure}
In Table \ref{CFD parameters} we can find the problem data that we used in our simulation. 
After we obtain a fully developed Karman vortex, and after the cylinder starts to rotate, we have a noticeable change in the direction of propagation of the vortex. This change of direction may hinder the representation of the solution by a small number of basis functions,
and this expectation is confirmed by running a POD on the fluid velocity snapshots collection, as we can see from Figure \ref{POD_Karman}.
\subsection{Preprocessing step}
We will focus on the preprocessing of the fluid velocity $\mathbf{u}_f$, since it is more straightforward to visualize the direction of propagation of the Karman vortex and hence understand the idea beneath the deformation map.
\begin{figure}
\centering
\includegraphics[scale=0.2]{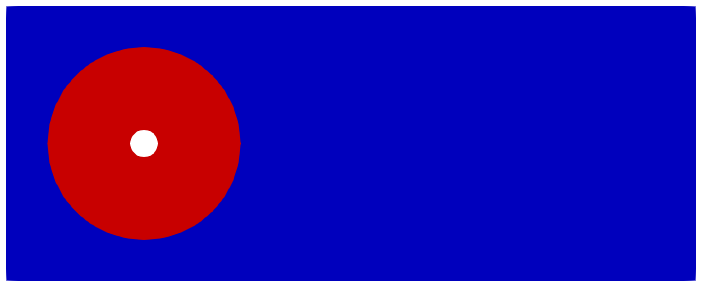}
\caption{Subdivision of the physical domain $\Omega$ into $\Omega_{int}$ (red) and $\Omega_{ext}$ (blue).}
\label{domain subdivision}
\end{figure}
First of all, let us notice from Figure \ref{domain subdivision} that, when building the mesh, we have defined a fictitious subdomain $\Omega_{int}$ (in red).
It is very important to choose the radius $r$ of the circular subdomain $\Omega_{int}$ in such a way that $\Omega_{int}$ is entirely contained in the physical domain $\Omega$ and yet $\Omega_{int}$ is able to capture all the
complex behaviour of the solution that is strictly related to the rotation of the cylinder. 
In fact the blue subdomain $\Omega_{ext}$ in Figure \ref{domain subdivision} deals with the small vortexes that originate in the wake of the cylinder before it starts to rotate: these vortexes propagate in the domain and their behaviour is not heavily
affected by the rotation of the cylinder, and therefore we should be able to reproduce their behaviour with a coarse set of basis functions.
Therefore as the solution features the most interesting phenomena in a neighborhood of the cylinder $\Gamma_{cyl}$, we are not going to preprocess the entire snapshot, 
but we are going to focus instead just on its restriction on $\Omega_{int}$. 
This procedure is more adapted, in the sense that we are not considering the solution manifold as a global entity, but we think of it as made of two manifolds, $\mathcal{M}_{ext}$ and $\mathcal{M}_{int}$: this division is based
on the division of the physical domain of interest $\Omega$ into two subdomains $\Omega_{ext}$ and $\Omega_{int}$. 
Therefore:
\begin{align}
\mathcal{M}_{ext} &= \{\mathbf{u}_f(\cdot, t)|_{\Omega_{ext}}, \quad t \in [0, T]\},\\
\mathcal{M}_{int} &= \{\mathbf{u}_f(\cdot, t)|_{\Omega_{int}}, \quad t \in [0, T]\}.
\end{align}
Being $\Omega_{ext}$ far from the rotating cylinder, and therefore far from the rotating phenomenon, we expect $\mathcal{M}_{ext}$ to be a better behaved solution manifold with respect to $\mathcal{M}_{int}$. This means that $\mathcal{M}_{ext}$
has a small Kolmogorov $n$-width, and does not need any preprocessing. On the contrary we will focus on $\mathcal{M}_{int}$, which will have a slowly decaying Kolmogorov $n$-width. Before going any further, let us remark that the 
subdivision of $\mathcal{M}_{\mathbf{u}_f}$ in $\mathcal{M}_{int}$ and $\mathcal{M}_{ext}$ is strictly dependent on the subdivision of $\Omega$, therefore the fictitious cylinder has to be chosen wisely, in such a way that $\Omega_{int}$ is able
to capture most of the rotating phenomenon. In particular, in our simulations, we have chosen $\Omega_{int}$ to be a cylinder of radius $r$, $7$ times larger than the radius of the physical cylinder.

For the preprocessing step therefore we first restrict the snapshots to the subdomain $\Omega_{int}$. Then, we want to build a (one-parameter) family of smooth and invertible maps
\begin{equation*} 
\mathcal{F}_{\mathbf{u}_f} = \{F_\theta\colon\Omega_{int}\to\Omega_{int},\ \theta = \theta(t),\ t \in [0, T]\},
\end{equation*}
where $\theta$ is the parameter identifying each map. 
In our case it is natural to choose, at each time $t$, $\theta(t)$ to be the angle spanned by the direction of propagation of the vortex  (obtained through a postprocessing of the solution $\mathbf{u}_f(t)$) and the horizontal axis. We choose the following preprocessing map:
\begin{equation}
\label{rotation}
F^{-1}_\theta(x, y) = 
\begin{pmatrix}
\cos\theta & -\sin\theta\\
\sin\theta & \cos\theta
\end{pmatrix}
\begin{pmatrix}
x \\
y
\end{pmatrix}
.
\end{equation}
The idea behind \eqref{rotation} is pretty straightforward: at each timestep $t$ of our simulation, we compute how much the vortex has changed its direction of propagation, and we therefore find $\theta(t)$. After that, in the
preprocessing step we take all the snapshots $\textbf{u}_f(t)$, we restrict them to the subdomain $\Omega_{int}$, and then we rotate them back to a horizontal direction of propagation of the vortex with $\textbf{u}_f(t) \circ F^{-1}_{\theta(t)}$.



\begin{figure}
\centering
\includegraphics[scale=0.5]{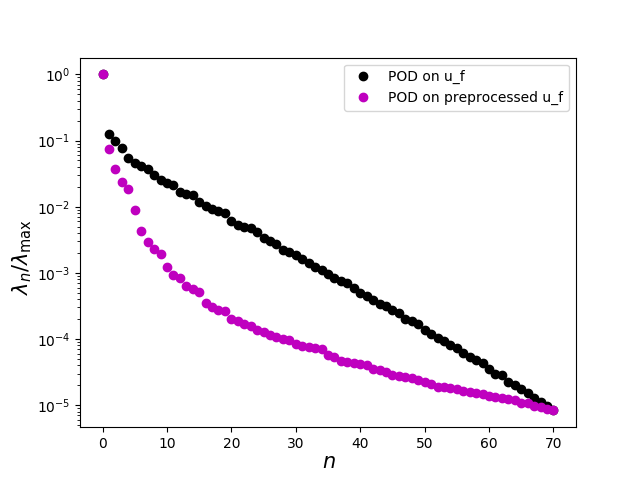}
\caption{POD comparison before (black) and after (magenta) the preprocessing of $\textbf{u}_f$.}
\label{POD comparison Karman vortex}
\end{figure}
In Figure \ref{POD comparison Karman vortex} we see the improvements (in terms of POD eigenvalues decay) that we get by applying the rotation to the snapshots: as we can see, after the preprocessing we obtain an improvement e.g. of almost $2$ orders of magnitude comparing the $10$-th eigenvalue of standard and preprocessed manifolds.

\section{A multiphysics problem}
We are now interested in applying the preprocessing procedure to a fluid-structure interaction problem, whose solution exhibits a transport dominated behaviour. The problem formulation features an Arbitrary Lagrangian Eulerian (ALE) formulation, see Richter \cite{Richter} for more details. An extensive explanation on how to treat such a coupled problem in a reduced order modelling setting can be found in Ballarin et al.\cite{BallarinRozza} (monolithic approach) and in Ballarin et al.\cite{BallarinRozzaMaday} (partitioned scheme).
For further details on reduced order models and applications of FSI problems we refer to Bertagna et al.\cite{BertagnaVeneziani}, Lassila et al. \cite{LassilaQuarteroniRozza,LassilaRozza} and Colciago \cite{Colciago}.
\subsection{Problem formulation}
We have a two dimensional rectangle of height $h_f$ and length $L$, filled with a Newtonian fluid. The structure is made by the top and the bottom of the rectangle, which are considered to 
be deformable, and their thickness is negligible with respect to the height of the rectangle. Since the structure is thin, it is described by a one dimensional model. 
We further assume that the displacement of the walls in the horizontal direction is negligible, and hence the structure presents only vertical motion: the behaviour of the compliant walls is therefore described by the
\emph{generalized string model} \cite{String_model_1,String_model_2}. 
We want to describe the behaviour of the solution (and the domain itself) in the time interval $[0, T]$.

\begin{figure}
\centering
\includegraphics[scale=0.2]{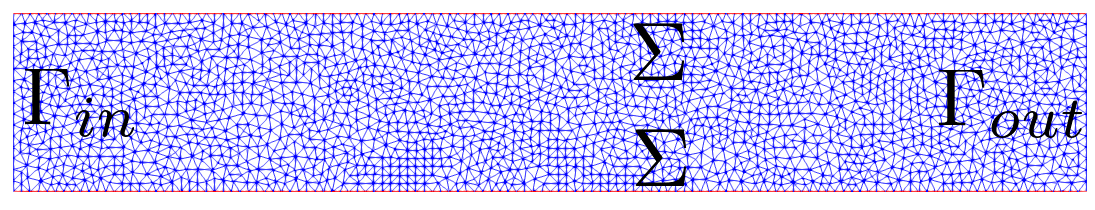}
\caption{Physical domain: fluid subdomain (blue) and structure subdomain (red). The fluid-structure interface coincides with the structure in our case.}
\label{fig:domain_FIS}
\end{figure}

Let $\Sigma$ denote the structure, and let $\Omega_t^f$ be the fluid domain at time $t$. The fluid equations are formulated, in this particular case, on a moving fluid domain $\Omega_t^f$. 
Let $\hat\Omega^f$ be the fluid reference domain: for convenience we take $\hat\Omega^f = \Omega^f = \Omega^f_{t=0}$ (the blue fluid domain in Figure \ref{fig:domain_FIS}).
\\The map $\mathbf{T}\colon\hat\Omega^f\to\Omega_t^f$ that maps $\hat\Omega^f$ to $\Omega_t^f$ at time $t$ is $\mathbf{T}$:
\begin{align*}
\mathbf{T}\colon\hat\Omega^f &\to\Omega_t^f\quad :\\
(x, y) &\mapsto (x, y+ d_f),
\end{align*}
where $d_f$ denotes the vertical component of the displacement (or deformation) of the fluid mesh, as such displacement in the horizontal direction is negligible, since the horizontal displacement of the structure is negligible.  
\\Thanks to the introduction of $\mathbf{T}$ we can therefore map the fluid equations back to the reference domain $\hat\Omega^f$, thus introducing additional terms to the classical Navier--Stokes equations; this results in an Arbitrary Lagrangian Eulerian (ALE)
formulation \cite{Richter, Gurtin} of our original problem.
Since we are considering the case of very small deformations, we will define $d_f$ as an harmonic extension of the structure displacement $d_s$:
\begin{equation}
\begin{cases}
-\Delta d_f = 0 &\text{in $\Omega^f$}, \\
d_f = d_s &\text{in $\Sigma$}.
\end{cases}
\label{eq:ext}
\end{equation}
There are alternative ways of defining the fluid displacement \cite{Richter} (pseudo-elasticity, bi-harmonic).

Let now $F$ be the gradient of $\mathbf{T}$, and let $J$ be the Jacobian.
The fluid-structure interaction problem reads: find fluid velocity $\textbf{u}_f(\cdot; t)\colon\Omega_t^f\to \mathbb{R}^2$, fluid pressure $p_f(\cdot; t)\colon\Omega_t^f\to\mathbb{R}$ 
and structure displacement $d_s(\cdot; t)\colon\Sigma\to\mathbb{R}$ such that:
\begin{equation}
\label{FSI_reference_domain}
\begin{cases}
J\rho_f(\partial_t\textbf{u}_f+F^{-1}(\textbf{u}_f-\partial_t d_f\textbf{e}_y)\cdot\nabla)\textbf{u}_f)-\text{div}(J\sigma^fF^{-T}) = \textbf{b$_f$} &\text{in $\Omega^f \times [0, T]$},\\
\text{div}(JF^{-1}\textbf{u}_f) = 0 &\text{in $\Omega^f\times[0, T]$},\\
\rho_s h_s\partial_{tt} d_s - c_0\partial_{xx} d_s + c_1 d_s = -\sigma^f(\textbf{u}_f, p_f)\textbf{n}\cdot\textbf{n} &\text{in $\Sigma$}.
\end{cases}
\end{equation}
Here $\rho_f$ is the fluid viscosity, $\rho_s$ is the structure viscosity, $h_s$ is the structure height (or thickness), $c_0$ and $c_1$ are the structure constitutive parameters. 
$\sigma^f$ is the Cauchy stress tensor for the fluid, and is defined as:
\begin{equation*}
\sigma_f(\textbf{u}_f, p_f) = -p_fI + \rho_f\nu_f\left(\nabla\textbf{u}_f F^{-1} + F^{-T}\nabla^T\textbf{u}_f\right).
\end{equation*} 
$\nu_f$ being the kinematic viscosity of the fluid and $I$ being the $2\times 2$ identity matrix. Finally, $\textbf{e}_y$ is the unit vector $(0, 1)$.  

Since it is a fluid-structure interaction problem, we need some coupling conditions. Let us denote $\Sigma(t)$ the fluid structure interface at time $t$; we have: 
\begin{equation}
\begin{cases}
\label{coupling conditions}
d_f=d_s &\text{in $\Sigma(t)$}, \qquad \text{continuity of the displacement},\\
\textbf{u}_f=\partial_t d_s\textbf{e}_y &\text{in $\Sigma(t)$}, \qquad\text{continuity of the velocity}.
\end{cases}
\end{equation}
The condition on the continuity of the displacement is a geometric condition, which stems from the fact that we do not want the fluid and the solid subdomain to overlap, and which we already employed in the definition of \eqref{eq:ext}.
As far as the second condition concerns, let us remark a very important fact: the quantity $\partial_t d_f$ does \emph{not} coincide with the fluid velocity $\textbf{u}_f$, 
as $\partial_t d_f$ denotes the \emph{velocity of the mesh}, and therefore it is not the physical velocity of the fluid particles.
For this reason, in general situations in FSI like our problem, we ask for a continuity condition for the physical fluid velocity $\textbf{u}_f$ and the mesh velocity $\partial_t d_f$ at the fluid--structure interface, 
that is exactly where the mesh velocity coincides with the structure velocity (see \eqref{eq:ext}$_2$).
Finally, there is a third coupling condition, which does not appear together with the others, since it is included in the right hand side of the generalized string equation \eqref{FSI_reference_domain}$_3$; 
this condition arises from the classical action-reaction principle. 

Together with the coupling conditions we also have to give some boundary and some initial conditions.
For the latter we assume that at time $t=0$ the system is at rest; the boundary conditions can be summarized in the following system:
\begin{equation}
\begin{cases}
\label{boundary conditions}
\sigma^f(\textbf{u}_f, p_f)\textbf{n} = -p_{in}(t)\textbf{n} &\text{in $\Gamma_{in}\times (0, T]$},\\
\sigma^f(\textbf{u}_f, p_f)\textbf{n} = -p_{out}(t)\textbf{n} &\text{in $\Gamma_{out}\times (0, T]$},\\
d_s = 0 &\text{in $\partial\Sigma\times [0, T]$,}
\end{cases}
\end{equation}
where $\textbf{n}$ is the outward normal. The third condition says that the structure is fixed at its extremities.\\
As far as the problem formulation concerns, let us remark that we adopted a supremizer enrichment technique also in this multiphysics test case, always to obtain, in the POD framework, a set of basis functions that allow for a stable approximation of the fluid pressure.
\subsection{Transport dominated FSI problem}
\begin{table}
\begin{center}
\begin{tabular}{|l|c||l|c|}
\hline
Data & Value & Data & Value\\
\hline
$\rho_f$ & $1$ g/cm$^3$ & $E_s$ & $0.75\times 10^6$ dyn/cm$^2$\\ 
\hline
$\nu_f$ & $0.035$ Poise & $\nu_s$ & $0.5$\\ 
\hline 
$b_f$ & $0$ & $c_0$ & $\frac{h_sE_s}{2(1+\nu_s)}$\\
\hline
$\rho_s$ & $1.1$ g/cm$^3$ & $c_1$ & $\frac{h_sE_s}{h_f^2(1-\nu_s^2)}$\\
\hline
$h_s$ & $0.1$ cm & $T_{in}$ & $2.5\times10^{-3}$\\
\hline
$\Delta T$ & $10^{-4}$ & $p_{in}(t)$ & $10^3\times[1-\text{cos}(\frac{2\pi t}{T_{in}})]\chi_{[0, 0.0025]}$\\
\hline
$N_{max}$ & $110$ & $p_{out}(t)$ & $\emptyset$\\
\hline
\end{tabular}
\end{center}
\caption{Problem data for the FSI test case}
\label{Constitutive parameters}
\end{table}
Problem data used for the simulation of our test case can be found in Table \ref{Constitutive parameters}: corresponding values are taken from the numerical results presented by Sy and Murea\cite{SyMurea_1,SyMurea_2}. 

\begin{figure}
\centering
\includegraphics[scale=0.35]{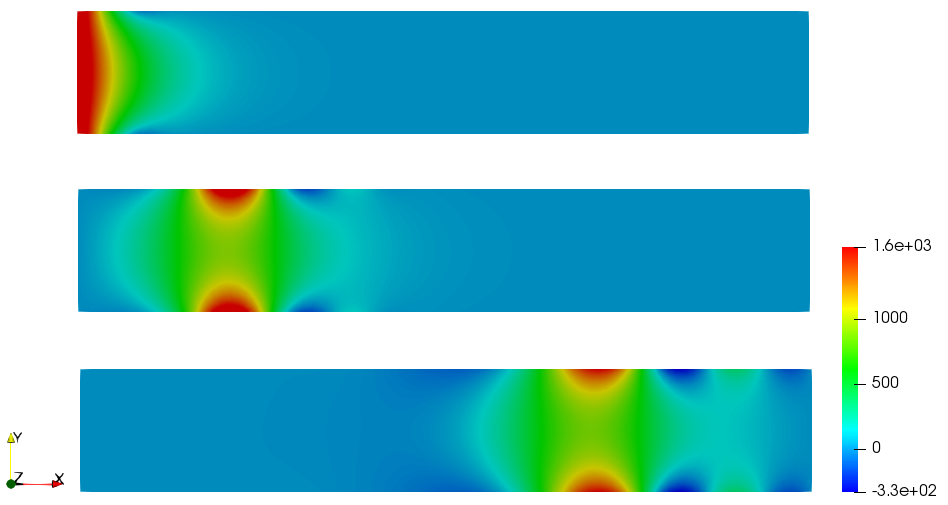}
\caption{Fluid pressure behaviour: the solution is pictured here at time $t=0.001$, $t=0.004$ and $t=0.011$. The peak of the wave is propagating into the domain, creating a transport phenomena.}
\label{fluid pressure behaviour}
\end{figure}

\begin{figure}
\centering
\includegraphics[scale=0.35]{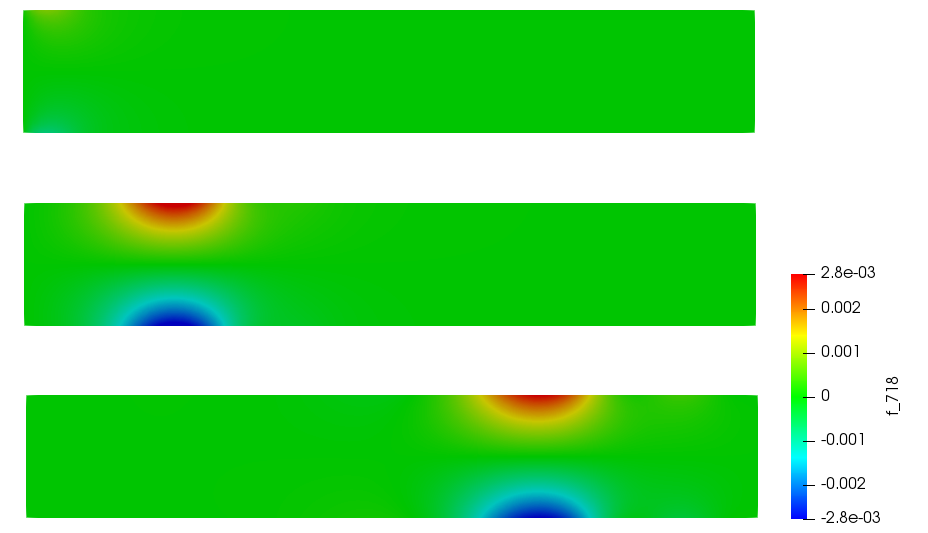}
\caption{Fluid displacement behaviour: again the solution is pictured at time $t=0.001$, $t=0.004$ and $t=0.011$. The peak of the wave is still very small at the beginning, it grows for some time and then it starts to propagate.}
\label{fluid displacement behaviour}
\end{figure}

\begin{figure}
\centering
\includegraphics[scale=0.35]{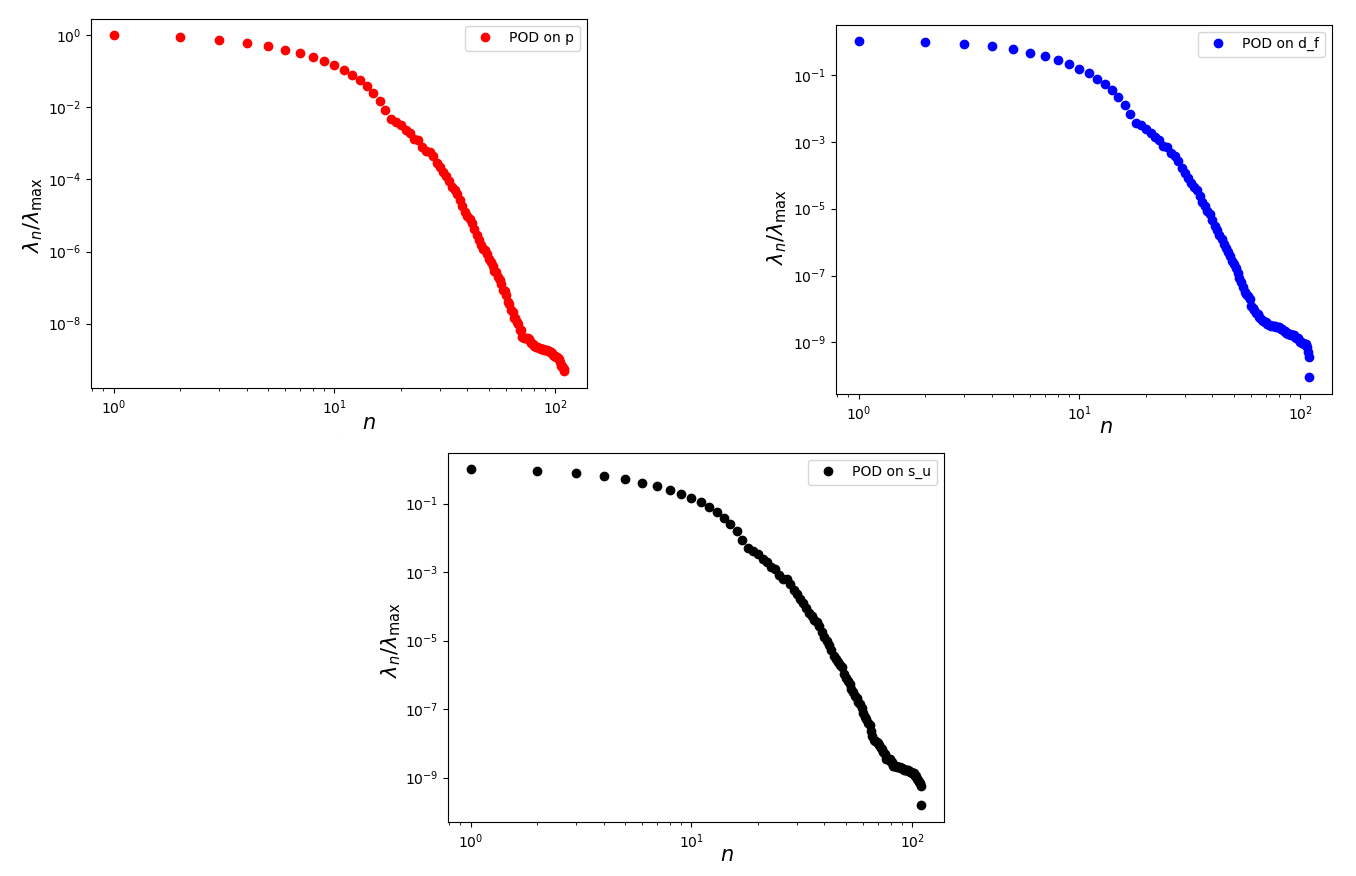}
\caption{Decay of the singular values for the POD on the fluid pressure (red, top left), fluid displacement (blue, top right) and fluid supremizer (black, bottom center).}
\label{POD_fsi}
\end{figure}

The behaviour of the fluid pressure $p_f$ and of the extended displacement $d_f$ is shown in Figure \ref{fluid pressure behaviour} and Figure \ref{fluid displacement behaviour} respectively. 
Our problem is transport dominated: if we look at Figure \ref{fluid pressure behaviour} for example, the change in time of the position of the peak of the pressure wave will be a difficult feature to capture at the reduced
order level with just a few modes. 
This expectation is finally confirmed at the numerical level, as we can see in Figure \ref{POD_fsi}, by running a POD on $p_f$, $d_f$, and also $s_u$.
Therefore we rely on a transformation on the set of solutions, in order to compensate this transport phenomenon.

In order to make the following exposition more clear and easy to read, from now on we focus only on a particular component of the solution of our problem, namely the fluid pressure. 
It is anyway important to keep in mind that, based on our simulation, all the components of the solution to the FSI problem are subject to a transport phenomenon, and hence every consideration that we are going to make on $p_f$
could be easily applied to any of the other components of the solution.

\subsection{Preprocessing step}
Let us see more in detail how to apply this preprocessing procedure to the fluid pressure solution manifold $\mathcal{M}_{p_f}$. Figure \ref{fluid pressure behaviour} shows how the peak of the pressure 
is transported in the domain. We would like to align the peaks at all time steps in a reference configuration; in this way, we obtain a set of snapshots where the 
pressure wave is not moving at all. In this case therefore a low number of modes will be sufficient to give a good representation of the situation. 

Starting from this observation we build a one parameter family of mappings 
\begin{equation*}
\mathcal{F}_{p_f} = \{ F_{\gamma}:\Omega^f \to \Omega^f, \quad \gamma=\gamma(t), \quad t\in [0, T]\}
\end{equation*}
such that for every $t$ in $[0, T]$, the peak of $p_f(F^{-1}_{\gamma}(\cdot), t)$ is located not at its original position, but is instead moved to the middle of the domain. 
In this way the new snapshots
$p_f(F^{-1}_{\gamma}(\cdot), t)$ will all have the peak located at the exact same position, meaning the middle of the domain.

When building the map $F^{-1}_{\gamma}$, another aspect to which we should pay great attention is the boundary conditions. 
Since we are not working in a periodic setting, we want to make sure that the preprocessed snapshots satisfy the same boundary conditions as the original snapshots. An easy way to make sure that these requirement is satisfied is to keep the points in $\Gamma_{in}$ and the points in $\Gamma_{out}$ fixed.
After some computations it follows that one possible map $F^{-1}_{\gamma}$ is the following:
\begin{equation*}
F^{-1}_{\gamma}(x) = \frac{3x\gamma}{x(\gamma-3)+3(6-\gamma)},
\end{equation*}
where $\gamma = \gamma(t)$ is the abscissa of the position of the peak of the wave at time $t$. We assume that the abscissa of the points on the inlet boundary $\Gamma_{in}$ is $x=0$ and the length of the domain $\Omega_f$ is $L=6$; 
in addition, the position in which we are moving the peak of the wave at every time $t$ is exactly in the center of the domain. Let us remark that the map $F^{-1}_t$ is just a stretching in the horizontal direction: 
this is due to the fact that we do not have any transport phenomena in the vertical direction, and hence there is no need for a transformation in the $y-$axis. So, with an abuse of notation, we can think of 
$F^{-1}_{\gamma}(x, y)$ as $F^{-1}_{\gamma}(x, y) = (F^{-1}_{\gamma}(x), y)$.

\begin{figure}
\centering
\includegraphics[scale=0.35]{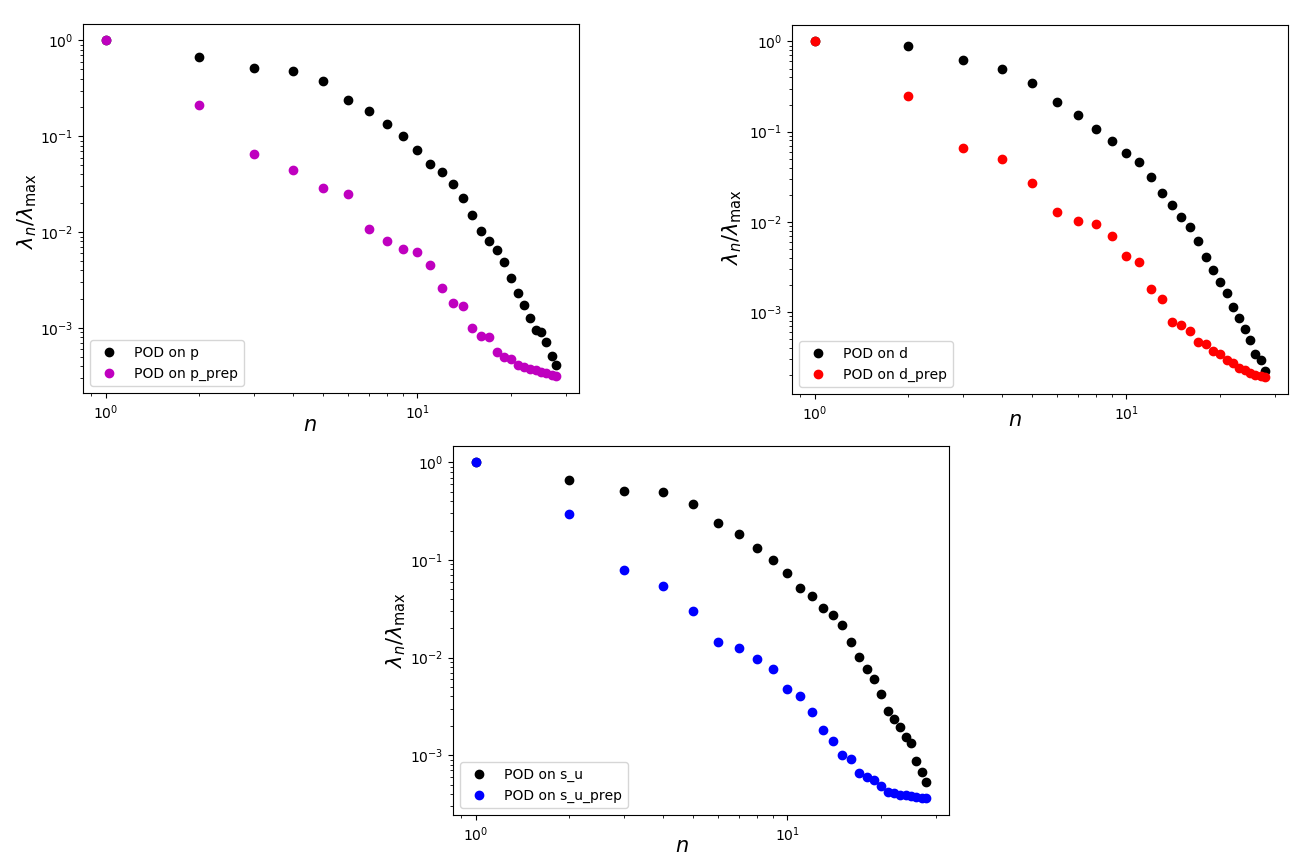}
\caption{Comparison on the decay of the eigenvalues for the POD on the original manifolds and on the preprocessed manifolds for the solutions components $p_f$ (top left), $d_f$ (top right) and $s_u$ (bottom center).}
\label{POD displacement improvements}
\end{figure}

\begin{figure}
\centering
\subfloat[][]{\includegraphics[width=0.51\textwidth]{fluid_pressure_behaviour.png}}\;
\subfloat[][]{\includegraphics[width=0.47\textwidth]{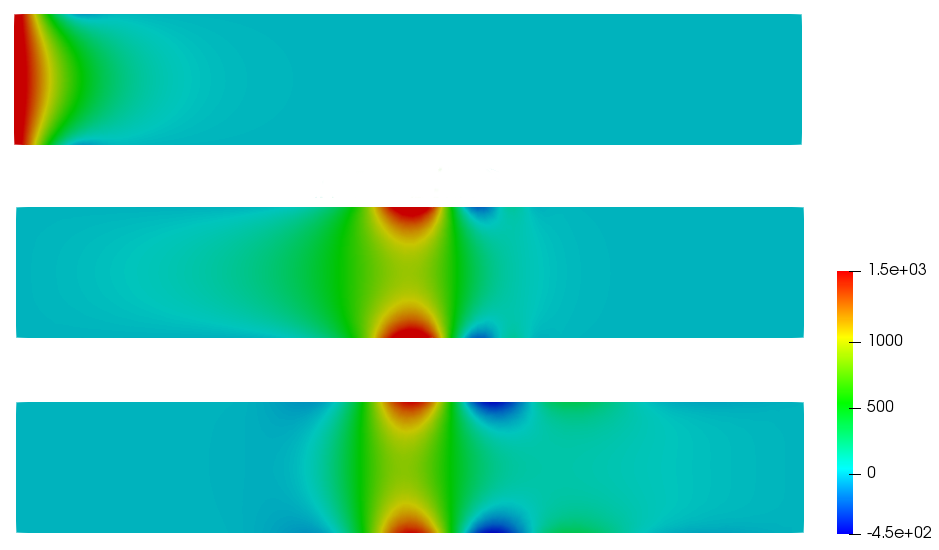}}
\caption{Original snapshots for $p_f$ at time $t=0.001$, $t=0.004$ and $t=0.011$ (a), and corresponding preprocessed snapshots (b).}
\label{snapshots_preprocessing}
\end{figure}
Our family of mappings $\mathcal{F}_{p_f}$ is a one-parameter family, therefore to identify the stretching map $F^{-1}_{\gamma}$ we have to compute the parameter $\gamma(t)$, which depends on time $t$.
First of all we compute all the snapshots $p_1, \dots, p_{N_{max}}$, where $N_{max} = \frac{T}{\Delta t}$, $p_i = p_f(t_i)$ and $t_i = i\Delta t$. Once we have the snapshots, 
we locate the abscissa $\gamma_i$ of the peak of the wave $p_i$; we then preprocess $p_i$ with $F^{-1}_{\gamma_i}$ and obtain a new snapshot $p_i(F^{-1}_{\gamma_i}(\cdot))$; see Figure \ref{snapshots_preprocessing} for a comparison between the original snapshots, 
and the preprocessed ones. 
The snapshot at time $t=0.001$ has not been preprocessed, and the reason for this is as follows: recall the expression of $p_{in}$ from Table \ref{Constitutive parameters}; as we can see, $p_{in}$ represents a given pulse at the inlet boundary, 
which is nonzero up to time $t=T_{in} = 0.0025$. So, up to time $t=0.0025$ the pulse makes the wave grow; when the pulse is null, the wave starts to propagate into the domain.
As previously mentioned, we are looking for a stretching map that moves the peak of the wave in a particular location of our choice, and 
at the same time keeps the points on the inlet boundary $\Gamma_{in}$ and on the outlet boundary $\Gamma_{out}$ fixed. However, in our simulation the peak of the wave is located exactly or very close to the inlet boundary for the
first $t=0.0024$ seconds: the transport phenomena is not immediate; this implies that, until time $t=0.0024$ we do not need to apply any preprocessing step. At the reduced level this translates into the fact that, up to $t=0.0024$
we adopt a standard reduction technique, using a standard reduced basis obtained with a standard POD on the first $24$ snapshots. From the snapshot $p_{25}$ to the last one, on the other hand, we first perform the preprocessing
step, and then perform a POD.

In Figure \ref{POD displacement improvements} we can compare the decay of the eigenvalues for the POD on the pressure with and without this preprocessing procedure. As we can see, with the preprocessing technique we do actually get an 
improvement in the decay of the eigenvalues, and in fact with less than $15$ modes we reach a level of $10^{-3}$, which is one order of magnitude less than the one we get with less than $15$ modes in the standard case.
We get the same results for the extended displacement $d_f$.

\section{CFD test case with a physical parameter}
So far we have considered test cases where the only parameter was time, 
nevertheless we are interested in investigating what happens if we add another physical parameter to the problem. 
We then go back to the CFD test case presented in section $3$, but now the Reynolds number $Re$ will be considered as a parameter. The fluid velocity solution manifold will be defined as:
\begin{equation*}
\mathcal{M}_{\textbf{u}_f} = \{\textbf{u}_f(t, \mu), \quad t\in[0, T], \quad \mu\in [Re_{\text{min}}, Re_{\text{max}}]\}, 
\end{equation*}
where now $\mu=Re$. Of course in this case we have to pay attention not only to the Reynolds number, but also to the rotation rate $\alpha$ of the cylinder, since it is strictly related to the development of a vortex shedding phenomenon.
For our particular CFD test case we choose the following parameter range:
\begin{equation*}
[Re_{\text{min}}, Re_{\text{max}}] = [47, 150],
\end{equation*}
and we choose $\alpha = 1.0$. We discretize the parameter space, choosing a set of parameter samplings $\{\mu_1,\dots,\mu_N\}$. For our problem we choose a Lagrange distribution sampling, i.e:
\begin{equation*}
\mu_i = Re_{\text{min}}\exp\left(\left(\frac{i-1}{N-1}\right)\log\left(\frac{Re_{\text{max}}}{Re_{\text{min}}}\right)\right).
\end{equation*}
After we have obtained a set of snapshots for each parameter in the parameter sampling set, we observe that
the change in the Reynolds number leads to changes in the behaviour of the fluid velocity, with the vortex shedding that may occur earlier or later, but in all the situations we see that after a while, due to the rotation
of the cylinder, the direction of propagation of the vortex changes. 
\subsection{POD-Greedy}
With the addition of a physical parameter, the exploration strategy will be carried out in a different way with respect to a standard POD on the set of snapshots; we are going to explore the parameter space with a pseudo-Greedy algorithm,
and we are going to explore in time with a POD on the set of snapshots corresponding to each parameter selected by the pseudo-Greedy strategy.
First of all we discretize the parameter space in order to obtain a sampling set of cardinality $N$ of our choice; in order to do so, we choose a Lagrange distribution sampling, i.e:
\begin{equation*}
\mu_i = Re_{min}\exp\left(\left(\frac{i-1}{N-1}\right)\log\left(\frac{Re_{max}}{Re_{min}}\right)\right).
\end{equation*}
Once we have the parameter sampling set $\{\mu_1, \dots, \mu_N\}$, we compute the truth solution for each one of these parameters. Afterwards, the POD is applied in the following way:
\begin{enumerate}
\item for $\mu_1$, we run a standard POD on the corresponding snapshots;
\item we now have at hand a set of reduced basis $\{\Phi_1,\dots,\Phi_{M_1}\}$;
\item for $\mu_i$, $i\geq 2$, we orthogonalize each snapshot $\textbf{u}_f^j(\mu_i)$ with respect to the linear space $\text{span}(\Phi_1, \dots, \Phi_{M_{i-1}})$;
\item we then run a POD on the set of orthogonalized snapshots, and add the resulting basis functions to the already existing set of basis functions.
\end{enumerate}
We remark that we are not actually using a proper Greedy algorithm, because we are choosing a priori the parameter sampling set; the reason for this ``pseudo''-Greedy is that it is beyond the scope of this work to focus on an error
estimator to be used for the Greedy algorithm, and also a simple Lagrange distribution will be sufficient to have an insight on what is going on before and after the preprocessing step. We also remark that different possibilities
for the POD are possible: the one we are using here is based on an orthogonalization step, where we are getting rid of the superfluous information before pursuing a POD \cite{Haasdonk}. Another possibility would be to run first a standard
POD on the snapshots, for every parameter, and then, at the end, run another standard POD on the set of obtained reduced basis, again to get rid of the superfluous information \cite{Nguyen}.
We do not present here any result showing the rate of decay of the greatest eigenvalue returned by running a POD with orthogonalization on the set of snapshots corresponding to each parameter of the sampling set. The reason of this 
choice is the fact that, as we previously remarked, we did not actually use a Greedy procedure to select the parameters in the parameter space; in addition, a plot showing the decay of the first eigenvalues of each POD gives a good
idea of how well the orthogonalization procedure is doing. The use of a proper Greedy algorithm based on a reliable error estimator and the proof that the orthogonalization procedure in the POD produces good advantages are outside
the scope of this paper. 
\subsection{Preprocessing step}
The preprocessing procedure is carried out, for every value of the parameter $\mu$ in the discretized parameter space, in the same way it was carried out in the case with no physical parameter: $\forall k \in \{\mu_1, \dots, \mu_N\}$, and
for every timestep $t_i$ we compute $\theta_i^k$, which is the angle spanned by the direction of propagation of the vortex and the horizontal axis. Then the deformation map $F^{-1}_{\theta_i^k}(x, y)$ is defined as:
\begin{equation*}
F^{-1}_{\theta_i^k}(x, y) = 
\begin{pmatrix}
\cos(\theta_i^k) & -\sin(\theta_i^k)\\
\sin(\theta_i^k) & \cos(\theta_i^k)
\end{pmatrix}
\begin{pmatrix}
x\\
y
\end{pmatrix}
.
\end{equation*}
Also in this case the idea is to rotate back to a horizontal direction of propagation all the vortexes. We remark that the preprocessing procedure is carried out on the snapshots restricted to the domain $\Omega_{int}$, exactly as
we did for the problem in Section $3$.
\begin{figure}
\includegraphics[scale=0.35]{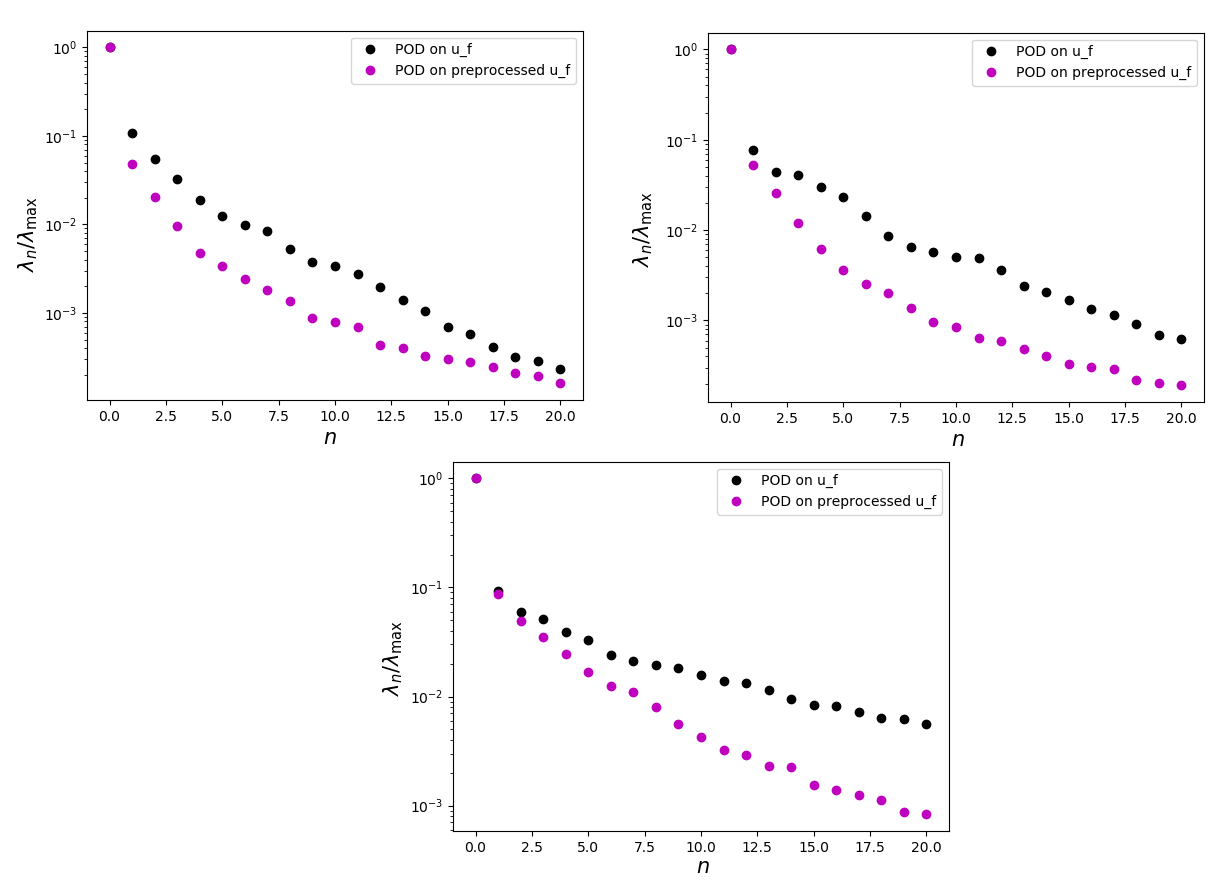}
\caption{Comparison of the rate of decay of the eigenvalues on the original velocity solution manifold $\mathcal{M}_{\textbf{u}_f}$ and on the preprocessed solution manifold $\mathcal{M}_{\textbf{u}_f}^p$ for $\mu=47$ (top left), $\mu=82.72$ (top right) and for $\mu=145$ (bottom).}
\label{parametric_results}
\end{figure}
Figure \ref{parametric_results} shows the results that we obtain for three different values of the Reynolds number. As we can see, there is an improvement in the rate of decay of the eigenvalues, with results showing a difference of one order of magnitude for $Re=145$ (right column) with just $20$ modes.

\section{Conclusions and future perspectives}
In our work we used a preprocessing of the snapshots during the offline stage of the reduced basis method to improve the rate of decay of the Kolmogorov $n$-width of the solution manifold of the problem of interest. 
We focused on two different test cases: a fluid dynamics problem
where a Karman vortex develops in the wake of a cylinder and where the vortex changes its direction of propagation due to the rotation of the cylinder, and a multiphysics problem, where the solution behaves like a travelling wave.
The preprocessing procedure is different according to the problem at hand, since the definition of the deformation maps changes according to the behaviour of the solution and to the boundary conditions used. 
The method performs very well for one dimensional problems \cite{tesidottoratocagniart}. In this work we focused on two dimensional problems in a non-periodic setting. 
We can say that by adopting this preprocessing procedure of the set of snapshots, the saves from a computational point of view in terms of the dimension of the set of basis functions needed to reach a certain 
approximation accuracy is evident.
The results that we obtained are promising: in the coupled problem, to reach a magnitude of $10^{-3}$ for the eigenvalues related to the POD on the 
pressure, we need $10$ less modes with respect to the standard situation with no preprocessing, and we can say the same about other components (fluid displacement, fluid velocity), thus lowering the dimension of the system to be solved in the online phase of the reduced method of at least $30$.
We obtained promising results also in the fluid dynamics test case. In the non parametric problem, results for the fluid velocity show that the eigenvalues after the preprocessing decay with almost two orders of magnitude faster than 
the standard case: to reach a magnitude of $10^{-3}$ in the standard case we need almost $30$ modes whereas in the preprocessed case we need $10$ modes. In the parametric case we analyzed the results for three different values of the 
physical parameter: all the three cases showed good results.
\\For future perspectives, we do believe that there are several possibilities to further develop and improve this preprocessing technique. 
In our work, both in the FSI problem and in the CFD one, we relied on our a priori knowledge of the behaviour of the solution in order to compute exactly the parameter $\gamma=\gamma(t)$ identifying the deformation map $F_{\gamma}$ at time $t$:
we computed exactly the position of the peak of the wave in the multiphysics problem, and we computed exactly the angle spanned by the direction of propagation of the vortex and the horizontal
axis in the fluid dynamic problem. There is thus a large space to develop methods to find a good parameter $\hat\gamma(t)$, instead of using the exact one, that ensure a good performance of the preprocessing technique and that are computationally
feasible.
In this work we did not go in the details of the online phase of the reduced method with the preprocessing technique, as we were mainly interested in studying the performance of the new method in different situations, and 
comparing the rate of decay of the eigenvalues to quantify the improvements that we obtain. Future research work will include: the efficient evaluation of the online phase, and, in addition, application of this reduced order method coupled with the preprocessing technique in the framework of inverse problems\cite{MadayMulaPateraYano, Galarce, StrazzulloBallarinRozza, StrazzulloBallarinMosettiRozza, Zainib}.

\subsection*{Acknowledgements}
We acknowledge the support by European Union Funding for Research and Innovation -- Horizon 2020 Program -- in the framework of European Research Council Executive Agency: Consolidator Grant H2020 ERC CoG 2015 AROMA-CFD project 681447 ``Advanced Reduced Order Methods with Applications in Computational Fluid Dynamics'' (PI Prof. Gianluigi Rozza).
We also acknowledge the INDAM-GNCS project ``Advanced intrusive and non-intrusive model order reduction techniques and applications'', $2019$.
The computations in this work have been performed with multiphenics \cite{multiphenics}, which is a library that aims at providing tools in FEniCS \cite{fenics} for an easy prototyping of multiphysics problems on conforming meshes, and RBniCS \cite{rbnics}, which is an implementation of several reduced order modelling techniques in FEniCS. We acknowledge developers and contributors to these three libraries.

\clearpage
\addcontentsline{toc}{section}{References}
\bibliographystyle{abbrv}
\bibliography{bibliography}
\end{document}